# ON TIME-INHOMOGENEOUS CONTROLLED DIFFUSION PROCESSES IN DOMAINS

By Hongjie Dong and N. V. Krylov[1]

*Institute for Advanced Study and University of Minnesota*

Time-inhomogeneous controlled diffusion processes in both cylindrical and noncylindrical domains are considered. Bellman's principle and its applications to proving the continuity of value functions are investigated.

The first part of this article is devoted to quite an old subject in the theory of controlled diffusion processes, namely deriving Bellman's principle (also called the principle of optimality) for processes controlled up to the first exit time from bounded domains. This principle plays a major role in many aspects of the theory of controlled diffusion processes. The necessity of proving it and deriving from it some continuity properties of value functions came to light while investigating the rate of convergence of finite-difference approximations for Bellman's equations. In this connection, we point out that our main results are Theorems 2.10, 2.13 and 2.17 in the second part of the paper. Theorem 2.17 is one of the main ingredients in [3], where we proved a sharp result that the rate of convergence of finite-difference approximations for Bellman's equations in bounded domains is not less than $h^{1/2}$, with $h$ being the mesh size. Theorem 2.17 is similar to Theorem 2.1 of [8] and is nontrivial even if we consider a single diffusion process without any control. In that case, it yields the rate of convergence $h^{1/2}$ without much work (see, e.g., Corollary 1.10 of [8]).

Our main results depend heavily on the validity of Bellman's principle. Bellman's principle has been derived in different settings in many papers and books. We refer the reader to [1, 2, 4, 5, 6, 10] and the references therein. Probably the article closest to the subject of the present one is [6], where

Received June 2005; revised November 2005.

[1]Supported in part by NSF Grant DMS-01-40405.

Bellman's principle is derived under very general conditions allowing unbounded domains and the coefficients of the controlled processes, but only for the problem of optimal stopping of controlled processes. Later on, these results were used to obtain sharp results concerning when the value functions for *time-homogeneous* processes satisfy the corresponding Bellman's equations. Another result, which is also very close to the results presented in the first part of this paper, is Theorem 2.1 in Chapter V of [4]. However, there are gaps in the original proof of the theorem (the corrected version is to appear in the forthcoming second edition of [4]) and the conditions under which it is stated are somewhat different from those we require for some applications we wish to consider. It is worth noting that in [4], the controlled process is considered up to the first exit time from the *closure* of a domain, so that if we have two domains with the same closure, the corresponding value functions will coincide. In contrast, we consider exit times from a domain as is usually done in the theory of Markov processes. One of major technical differences between these two settings is that our exit times are lower semicontinuous and the exit times from [4] are upper semicontinuous.

The approach in [4] originated from [11], where the reader can also find many useful results concerning the continuity of value functions.

In Section 1 we prove Bellman's principle in a setting more general than that of Theorem 2.1 in Chapter V of [4] (see Remark 1.14). Several examples show that under the assumptions in Section 1, value functions can be discontinuous even inside the domains. With additional assumptions, in Section 2, we prove the Lipschitz continuity of value functions in space variables and Hölder-1/2 continuity in the time variable, which is one of the main motivations of this paper. In Remark 2.14, we also present our understanding as to how the statement of Theorem 2.1 in Chapter V of [4] regarding the continuity of value functions can be corrected. Finally, we derive in Corollaries 1.3 and 2.12 an inequality which we use to prove Theorem 2.17. As we have mentioned above, the last theorem plays a major role in investigating the rate of convergence of numerical approximations for Bellman's equations in domains.

**1. Bellman's principle.** Let $A$ be a separable metric space and let $A(n)$ be fixed subsets of $A$, $n = 1, 2, \ldots$, such that $A = \bigcup_n A(n)$, $A(n) \subset A(n+1)$. Let $(\Omega, \mathcal{F}, P)$ be a complete probability space and $\{\mathcal{F}_t; t \geq 0\}$ an increasing filtration of $\sigma$-algebras $\mathcal{F}_t \subset \mathcal{F}$ which are complete with respect to $\mathcal{F}$, $P$. Let $(w_t, \mathcal{F}_t; t \geq 0)$ be a $d_1$-dimensional Wiener process on $(\Omega, \mathcal{F}, P)$.

Suppose that the following have been defined for $\alpha \in A$ and $(t, x) \in \mathbb{R} \times \mathbb{R}^d$: a $d \times d_1$ matrix $\sigma^\alpha(t, x)$, a $d$-dimensional vector $b^\alpha(t, x)$ and real numbers $c^\alpha(t, x)$, $f^\alpha(t, x)$ and $g(t, x)$. We assume that for every $n \geq 1$, on $A(n) \times \mathbb{R}^{d+1}$, the functions $\sigma, b, c$ and $f$ are Borel, bounded, continuous in $(\alpha, x)$ and continuous in $x$ uniformly with respect to $\alpha$ for each $t \in \mathbb{R}$. Moreover, for every $n \geq 1$, on $A(n) \times \mathbb{R}^{d+1}$, let $\sigma$ and $b$ satisfy a Lipschitz condition in $x$ with constant not depending on $(\alpha, t)$ and let $g$ be lower semicontinuous and bounded in $\mathbb{R}^{d+1}$.

By $\mathfrak{A}(n)$, we denote the set of all functions $\alpha_r(\omega)$ on $\Omega \times [0, \infty)$ which are $\mathcal{F}_r$-adapted and measurable in $(\omega, r)$ with values in $A(n)$. Let $\mathfrak{A} = \bigcup_n \mathfrak{A}(n)$ and let $\mathfrak{M}$ be the set of all bounded stopping times (relative to $\{\mathcal{F}_r\}$).

For $\alpha \in \mathfrak{A}$ and $(t, x) \in \mathbb{R}^{d+1}$, we consider the Itô equation

$$(1.1) \qquad x_s = x + \int_0^s \sigma^{\alpha_r}(t + r, x_r) \, dw_r + \int_0^s b^{\alpha_r}(t + r, x_r) \, dr.$$

The solution of this equation is known to exist and to be unique. We denote this solution by $x_s^{\alpha, t, x}$, following the abbreviated notation adopted in [5].

For any $s \geq 0$, we set

$$\varphi_s = \varphi_s^{\alpha, t, x} = \int_0^s c^{\alpha_r}(t + r, x_r^{\alpha, t, x}) \, dr.$$

Let $Q$ be a *bounded* domain in $\mathbb{R}^{d+1} = \mathbb{R} \times \mathbb{R}^d$ and let $\tau = \tau^{\alpha, t, x}$ be the first exit time of $(t + s, x_s^{\alpha, t, x})$ from $Q$:

$$\tau^{\alpha, t, x} = \inf\{s \geq 0 : (t + s, x_s^{\alpha, t, x}) \notin Q\}.$$

Observe that since $Q$ is bounded, $\tau^{\alpha, t, x}$ is a *bounded* stopping time.

Define the parabolic boundary $\partial' Q$ of $Q$ as the set of all points $(t, x)$ on $\partial Q$ for each of which there exists a curve $(s, y_s)$, $t - \varepsilon \leq s \leq t$, such that $\varepsilon > 0$, $(t, y_t) = (t, x)$, $y_s$ is a continuous function and $(s, y_s) \in Q$, $t - \varepsilon \leq s < t$. Obviously, if $(t, x) \in Q$, then at $s = \tau^{\alpha, t, x}$, the point $(t + s, x_s^{\alpha, t, x})$ lies on $\partial' Q$.

Set

$$v^\alpha(t, x) = E_{t,x}^\alpha \left[ \int_0^\tau f^{\alpha_s}(t + s, x_s) e^{-\varphi_s} \, ds + g(t + \tau, x_\tau) e^{-\varphi_\tau} \right],$$

$$v = \sup_{\alpha \in \mathfrak{A}} v^\alpha,$$

where we use common abbreviated notation, according to which we put the indices $\alpha$, $t$, $x$ beside the expectation sign instead of explicitly exhibiting them inside the expectation sign for every object that can carry all or part of them. For instance,

$$E_{t,x}^\alpha g(t + \tau, x_\tau) e^{-\varphi_\tau} = E g(t + \tau^{\alpha, t, x}, x_{\tau^{\alpha, t, x}}^{\alpha, t, x}) \exp(-\varphi_{\tau^{\alpha, t, x}}^{\alpha, t, x}).$$

It is worth noting that $\tau^{\alpha, t, x} = 0$ if $(t, x) \notin Q$. Therefore, $v = g$ in $Q^c$.

The above assumptions and notation will apply throughout the paper. Additional assumptions will be introduced for each particular result.

Observe that since $Q$ is bounded, $v \geq -N$, where $N$ is a constant, and the case $v \equiv \infty$ in $Q$ is not excluded.

The following version of Bellman's principle is the first main result of this section.

THEOREM 1.1. *Assume that $g \equiv 0$ and that there is an $\underline{\alpha} \in A$ such that $f^{\underline{\alpha}} \geq 0$ on $Q$. Then:*

(i) *the function $v$ is Borel measurable, nonnegative and, moreover, it is lower continuous in $Q$, that is,*

$$v(t,x) = \liminf_{(s,y) \to (t,x)} v(s,y) \qquad \forall\, (t,x) \in Q;$$

(ii) *we have*

$$(1.2) \quad v(t,x) = \sup_{\alpha \in \mathfrak{A}} E_{t,x}^{\alpha} \left[ \int_0^{\gamma} f^{\alpha_s}(t+s, x_s) e^{-\varphi_s}\, ds + v(t+\gamma, x_{\gamma}) e^{-\varphi_{\gamma}} \right]$$

*whenever $(t,x) \in \bar{Q}$ and for any $\alpha \in \mathfrak{A}$, we are given a stopping time $\gamma^{\alpha} \leq \tau^{\alpha,t,x}$ [in (1.2) the superscript $\alpha$ of $\gamma^{\alpha}$ is dropped in accordance with the above stipulation].*

PROOF. For any $\gamma \in \mathfrak{M}$, we set

$$v^{\alpha,\gamma}(t,x) = E_{t,x}^{\alpha} \int_0^{\gamma} f^{\alpha_s}(t+s, x_s) e^{-\varphi_s}\, ds,$$

$$w(t,x) = \sup_{\alpha \in \mathfrak{A}} \sup_{\gamma \in \mathfrak{M}} v^{\alpha, \gamma \wedge \tau}(t,x).$$

It is known from [6] (see Theorems 1.1, 2.4 and Lemma 2.2 therein) that under the conditions of the theorem, $w \geq 0$, the function $w$ is Borel, it is also lower continuous in $Q$, we have

$$(1.3) \qquad\qquad\qquad w(t,x) = v(t,x)$$

on $\bar{Q}$ and the process

$$\rho_s = \rho_s^{\alpha,t,x} = w(t+s \wedge \tau, x_{s \wedge \tau}) e^{-\varphi_{s \wedge \tau}} + \int_0^{s \wedge \tau} f^{\alpha_r}(t+r, x_r) e^{-\varphi_r}\, dr$$

is a supermartingale on $[0, \infty)$ for any $\alpha \in \mathfrak{A}$. Therefore, we have

$$(1.4) \qquad\qquad v(t,x) = E_{t,x}^{\alpha} \rho_0 \geq E_{t,x}^{\alpha} \rho_{\gamma} \geq E_{t,x}^{\alpha} \rho_{\tau}.$$

After taking supremum over $\alpha \in \mathfrak{A}$ in (1.4), we obtain

$$(1.5) \qquad v(t,x) = \sup_{\alpha \in \mathfrak{A}} E_{t,x}^{\alpha} \rho_0 \geq \sup_{\alpha \in \mathfrak{A}} E_{t,x}^{\alpha} \rho_{\gamma} \geq \sup_{\alpha \in \mathfrak{A}} E_{t,x}^{\alpha} \rho_{\tau}.$$

Since the rightmost term in (1.5) equals $v(t,x)$ by definition, all the inequalities in (1.5) are equalities. To complete the proof of (1.2), it only remains to observe that its right-hand side equals $\sup_{\alpha} E_{t,x}^{\alpha} \rho_{\gamma}$. $\quad\square$

REMARK 1.2. Since $v = 0$ on $\partial Q$ (even in $Q^c$) and $v \geq 0$ in $Q$, the lower continuity of $v$ holds on $\bar{Q}$, provided that $\partial Q$ has no isolated points, because $v$ is *continuous* along $\partial Q$, being identically zero there.

Set

$$D_t = \frac{\partial}{\partial t}, \qquad D_i = \frac{\partial}{\partial x^i}, \qquad D_{ij} = D_i D_j, \qquad a = (1/2)\sigma\sigma^*,$$

$$L^\alpha = L^\alpha(t,x) = (a^\alpha)^{ij}(t,x)D_{ij} + (b^\alpha)^i(t,x)D_i - c^\alpha(t,x).$$

As a corollary of Theorem 7.4 of [7] (or of the corresponding results in [5]) and Theorem 1.1, we have the following result. We remind the reader that Lipschitz continuous functions have bounded first-order generalized derivatives.

COROLLARY 1.3.   *In addition to the assumptions of Theorem* 1.1, *suppose that $v$ is bounded in an open set $Q' \subset Q$ and that for any $\alpha \in A$, the generalized function*

$$(1.6) \qquad\qquad \sum_{i,j} D_{ij}(a^\alpha)^{ij}$$

*is a locally integrable function on $Q'$. Then*

$$(1.7) \qquad\qquad D_t v + L^\alpha v + f^\alpha \le 0$$

*in $Q'$ in the sense of generalized functions, that is, for any nonnegative $\chi \in C_0^\infty(Q')$,*

$$\int_{Q'} v(-D_t + L^{\alpha*})\chi \, dt \, dx + \int_{Q'} f^\alpha \chi \, dt \, dx \le 0,$$

*where*

$$L^{\alpha*} := (a^\alpha)^{ij} D_{ij} - (b^\alpha)^i D_i + 2[D_j(a^\alpha)^{ij}]D_i + D_{ij}(a^\alpha)^{ij} - D_i(b^\alpha)^i - c^\alpha.$$

EXAMPLE 1.4.   Generally, in the situation of Theorem 1.1, the function $v$ need not be continuous in $Q$, even if $f$ is bounded. For instance, take $d = 2$ and consider the following (uncontrolled deterministic) process in $\mathbb{R}^2 = \{(x,y) : x, y \in \mathbb{R}\}$:

$$(1.8) \qquad\qquad dx_t = dt, \qquad dy_t = 0.$$

Let $c = 0$, $f = 1$ and $Q = (-1,4) \times (B_2 \setminus \bar{B}_1)$, where $B_r$ is the open ball of radius $r$ centered at the origin. As is easily seen, $v(0,x,y)$ is discontinuous along the lines $y = \pm 1$, $x \in [-\sqrt{3}, 0]$.

To obtain a generalization of Theorem 1.1 for $g \not\equiv 0$, we need one more assumption on $Q$.

Assumption 1.5. There exists a function $\psi \in C(\bar{Q})$ such that the first derivatives of $\psi$ with respect to $(t, x)$ and the second derivatives with respect to $x$ are continuous on $\bar{Q}$, $\psi$ vanishes on the parabolic boundary $\partial' Q$ of $Q$ and for some $\underline{\alpha} \in A$,

$$(1.9) \qquad D_t \psi + L^{\underline{\alpha}} \psi \leq -1 \qquad \text{in } Q.$$

Theorem 1.6. *Suppose that Assumption* 1.5 *is satisfied. Then assertions* (i) *and* (ii) *of Theorem* 1.1 *hold true. In particular, if $v$ is bounded in an open set $Q' \subset Q$ and for any $\alpha \in A$, the generalized function* (1.6) *is a locally integrable function on $Q'$, then* (1.7) *holds in $Q'$ in the sense of generalized functions, as in Corollary* 1.3.

Proof. To exhibit the dependence of $v$ and $v^\alpha$ on $g$, write $v = v[g]$ and $v^\alpha[g]$. Since $g$ is bounded and lower semicontinuous, there exists a sequence of smooth functions $g_n \uparrow g$. Also, notice that by the monotone convergence theorem,

$$v^\alpha[g] = \sup_n v^\alpha[g_n],$$

implying that

$$(1.10) \qquad v[g] = \sup_\alpha \sup_n v^\alpha[g_n] = \sup_n \sup_\alpha v^\alpha[g_n] = \sup_n v[g_n].$$

This shows how to obtain (1.2) for $g$ from the same assertion for $g_n$, so that in the proof of (1.2), we may assume that $g$ is a smooth function. Furthermore, since $\tau = \tau^{\alpha,t,x} = 0$ if $(t, x) \in \partial Q$, we may assume that $(t, x) \in Q$.

Next, let $N$ be a positive real number to be chosen later. Owing to Itô's formula and the facts that $(t + \tau, x_\tau^{\alpha,t,x}) \in \partial' Q$ and $\psi = 0$ on $\partial' Q$, we have

$$v^\alpha(t, x) - g(t, x) + N\psi(t, x)$$

$$(1.11) \qquad = E_{t,x}^\alpha \left[ \int_0^\tau (f^{\alpha_s} + D_t g + L^{\alpha_s} g \right.$$

$$\left. - N D_t \psi - N L^{\alpha_s} \psi)(t + s, x_s) e^{-\varphi_s} \, ds \right].$$

Denote by $\bar{f}^\alpha$ any continuous continuation of

$$D_t g + L^\alpha g - N D_t \psi - N L^\alpha \psi$$

outside $\bar{Q}$. This is possible because by assumption, the derivatives of $g$ and $\psi$ involved above are continuous in $\bar{Q}$. Then $v^\alpha - g + N\psi$ is simply $v^\alpha$ constructed from $f^\alpha + \bar{f}^\alpha$ and $g = 0$.

If (1.2) holds with

$$v - g + N\psi, \qquad f^{\alpha_s} + \bar{f}^{\alpha_s}$$

in place of $v, f^{\alpha_s}$, respectively, then by using Itô's formula again, we obtain (1.2) in its original form. This enables us to assume that $g = 0$. Due to (1.9), we can choose $N$ sufficiently large so that

$$f^{\underline{\alpha}} + \bar{f}^{\underline{\alpha}} = f^{\underline{\alpha}} + D_t g + L^{\underline{\alpha}} g - N D_t \psi - N L^{\underline{\alpha}} \psi \geq 0$$

in $Q$. Then (1.2) follows immediately from Theorem 1.1. This theorem also shows that $v$ is lower continuous in $Q$, at least for smooth $g$. Then (1.10) implies that $v$ is lower semicontinuous in $Q$ in the general case. The fact that it is lower continuous follows from (1.2), as in the proof of Theorem 2.4 of [6]. The theorem is thus proved. □

Next, we prove a similar result for processes in cylindrical domains. Let $D$ be a bounded domain in $\mathbb{R}^d$ and $T \in (0, \infty)$ a fixed number. Usually, one is interested in processes $(t + s, x_s^{\alpha,t,x})$ not until $\tau^{\alpha,t,x}$, but rather

$$(1.12) \qquad (T - t) \wedge \inf\{s \geq 0, x_s^{\alpha,t,x} \notin D\},$$

which is the first exit time of $(t + s, x_s^{\alpha,t,x})$ from $(-\infty, T) \times D$. These two exit times coincide if we take $(0, T) \times D$ as $Q$ and $t > 0$. However, if $t = 0$, then the former exit time is zero, since the starting point is already outside $(0, T) \times D$. Psychologically, the value $t = 0$ *looks* important and, therefore, in order to allow the process $(t + s, x_s^{\alpha,t,x})$ to start at points $(0, x)$ and yet have nontrivial objects to deal with, we set

$$Q = (-1, T) \times D.$$

We impose an assumption slightly different from Assumption 1.5:

Assumption 1.7. There exists a function $\psi \in C(\bar{Q})$ such that the first derivatives of $\psi$ with respect to $(t, x)$ and the second derivatives with respect to $x$ are continuous on $\bar{Q}$, $\psi > 0$ in $Q$ and $\psi$ vanishes on $(-1, T) \times \partial D$. Condition (1.9) is also satisfied for an $\underline{\alpha} \in A$.

Observe that in Assumption 1.7, we do not require $\psi$ to vanish on the whole parabolic boundary of $Q$. The reason is that if the derivatives of $\psi$ are continuous at points on $\{T\} \times \partial D$ and $\psi = 0$ on $\{T\} \times D$ and $(0, T) \times \partial D$, then the left-hand side of (1.9) is zero on $\{T\} \times D$ and so this inequality cannot be satisfied.

Remark 1.8. If in Assumption 1.7, inequality (1.9) holds only in $(-1, T) \times (D \setminus D')$, where $D' \subset \bar{D}' \subset D$, then one can modify $\psi$ in such a way that (1.9) holds in $Q$ for the modification. To see this, it suffices to observe that

$$(D_t + L^\alpha(t, x))[\psi(t, x) e^{\lambda(T - t)}] = e^{\lambda(T - t)}(D_t + L^\alpha)\psi(t, x) - e^{\lambda(T - t)}\lambda\psi(t, x)$$

and to choose $\lambda$ large enough, which would work if $\psi > 0$ in $[-1, T] \times \bar{D}'$.

It is also worth noting that Assumption 1.7 does not imply that $\partial D$ is smooth. For instance, if $D \subset \mathbb{R}^2 = \{(x, y) : x, y \in \mathbb{R}\}$ near the origin is described by $y > 2|x|$ and $L^\alpha \equiv \Delta$, then near the origin, one can take $\psi = y^2 - 4x^2$.

THEOREM 1.9. *Suppose that Assumption* 1.7 *is satisfied. Then assertions* (i) *and* (ii) *of Theorem* 1.1 *hold true. In particular, if $v$ is bounded in an open set $Q' \subset Q$ and for any $\alpha \in A$, the generalized function* (1.6) *is a locally integrable function on $Q'$, then* (1.7) *holds in $Q'$ in the sense of generalized functions, as in Corollary* 1.3.

PROOF. As in the proof of Theorem 1.6, we reduce the general situation to the case where $g$ is smooth and then, using Itô's formula, to the case $g = 0$. Let $\eta \in C^\infty(\mathbb{R})$ be a function satisfying

$$0 \le \eta \le 1 \qquad \text{in } \mathbb{R}, \qquad \eta \equiv 1 \qquad \text{in } (-\infty, -2], \qquad \eta \equiv 0 \qquad \text{in } [-1, +\infty).$$

For any $\varepsilon > 0$ and $\alpha \in A$, set

$$f_\varepsilon^\alpha(t, x) = f^\alpha(t, x)\eta(\varepsilon^{-1}(t - T))$$

and on $\bar{Q}$, define $v_\varepsilon^\alpha$ and $v_\varepsilon$ with $f_\varepsilon^\alpha$ in place of $f^\alpha$. From the definition of $v$, it is easy to see that $v_\varepsilon \to v$ uniformly on $\bar{Q}$ as $\varepsilon \downarrow 0$. Therefore, it suffices to prove the theorem for functions $f$ satisfying the additional assumption

$$f^\alpha \equiv 0 \qquad \text{on } [T - \varepsilon, T] \times \bar{D}$$

for some $\varepsilon > 0$ and any $\alpha \in A$.

Our goal is to apply Theorem 1.1. Denote

$$\bar{\psi} = (T - t)\psi/\varepsilon.$$

Obviously, $\bar{\psi} > 0$ in $Q$ and it vanishes on $\partial' Q$. In $(-1, T - \varepsilon) \times D$, we have

$$D_t \bar{\psi} + L^{\underline{\alpha}} \bar{\psi} \le \varepsilon^{-1}(T - t)(D_t \psi + L^{\underline{\alpha}} \psi) \le -1.$$

Meanwhile, in $[T - \varepsilon, T] \times D$, we have

$$D_t \bar{\psi} + L^{\underline{\alpha}} \bar{\psi} \le 0.$$

Therefore, we can choose $N$ sufficiently large such that in $Q$,

$$f^{\underline{\alpha}} - N D_t \psi - N L^{\underline{\alpha}} \psi \ge 0.$$

By using Theorem 1.1 and the argument in the proof of Theorem 1.6, we complete the proof of the present theorem. $\quad\square$

REMARK 1.10. It is known (see, e.g., [5]) that the optimal stopping problem for controlled diffusion processes reduces to a problem without stopping, but with the data $c^\alpha, f^\alpha$ becoming unbounded in the variable $\alpha$. Then the above results become applicable to the optimal stopping problem for controlled diffusion processes. This shows the usefulness of allowing our data to be unbounded in $\alpha$.

In the following example we present a situation in which $b^\alpha$ and $f^\alpha$ are unbounded:

EXAMPLE 1.11. Consider the so-called *singular stochastic control problem*. Let $x_t^\alpha$ be a process in $\mathbb{R}^d$ defined by

$$(1.13) \qquad x_t^\alpha = x + w_t + \alpha_t,$$

where $w_t$ is a $d$-dimensional Wiener process and $\alpha = \alpha_t$ is a $d$-dimensional control process such that for any $t \geq 0$, $\alpha_t$ is $\mathcal{F}_t$-measurable. Moreover, we will allow *any such continuous process* for which

$$|\alpha|_t := \mathrm{Var}_{[0,t]}\alpha$$

$$:= \sup\left\{\sum_{i=1}^{n-1} |\alpha_{t_{i+1}} - \alpha_{t_i}| : n = 1, 2, \ldots, 0 \leq t_1 < \cdots < t_n \leq t\right\}$$

$$< \infty \qquad \forall t,$$

that is, we allow processes of locally bounded total variation. Fix a smooth bounded domain $D \subset \mathbb{R}^d$, a lower semicontinuous bounded function $g = g(t, x)$ on $\mathbb{R} \times \mathbb{R}^d$ and a bounded continuous function $f$ on $\mathbb{R}^d$.

Assume that for $t \leq T$, we need to investigate

$$v(t, x) = \sup_\alpha E\left\{e^{-\tau}g(t + \tau, x_\tau) + \int_0^\tau e^{-s}f(x_s)\,ds - \int_0^\tau e^{-s}\,d|\alpha|_s\right\}$$

where $\tau$ is the minimum of the first exit time of $x_t$ from $D$ and $T - t$.

The fact that we restrict ourselves to continuous $\alpha_t$ allows us to use smooth approximations of $\alpha$ and to do this in such a way that the exit points for the original process and its approximations are close. Obviously, one can approximate process (1.13) by processes of the form

$$(1.14) \qquad x_t^\beta = x + w_t + \int_0^t \beta_s\,ds,$$

where $\beta \in \mathfrak{A} = \bigcup_n \mathfrak{A}(n)$ and $\mathfrak{A}(n)$ is defined as the set of jointly measurable $\mathcal{F}_t$-adapted processes with values in $A(n) = \{\beta \in \mathbb{R}^d : |\beta| \leq n\}$. It is also clear that

$$v(t, x) = \sup_{\beta \in \mathfrak{A}} E_x^\beta\left\{e^{-\tau}g(t + \tau, x_\tau) + \int_0^\tau e^{-s}(f(x_s) - |\beta_s|)\,ds\right\}.$$

Observe that due to the boundedness and smoothness of $D$, there is a smooth function $\psi = \psi(x)$ such that $\Delta \psi = -2$ in $D$ and $\psi = 0$ on $\partial D$. It follows that Assumption 1.7 is satisfied with $\underline{\alpha} = 0$. Therefore, Bellman's principle is applicable in this situation.

From Theorem 1.9, one can extract more information. Indeed, we have

$$D_t v + (1/2)\Delta v + \beta^i D_i v - v + f - |\beta| \leq 0$$

in $Q = (-1, T) \times D$ for any $\beta \in \mathbb{R}^d$. By considering large $|\beta|$ we see that in the sense of generalized functions, we have $\xi^i D_i v \leq 1$ in $Q$ for any unit $\xi \in \mathbb{R}^d$.

It is known that if the generalized gradient with respect to $x$ of a function that is measurable in $(t, x)$ is bounded by 1, then the function itself coincides [$(t,x)$-a.e.] with a function whose Lipschitz constant with respect to $x$ is majorated by 1. It follows that in $Q$ (a.e.), the function $v$ coincides with a function $\bar{v}$ that is Lipschitz continuous in $x$ with the Lipschitz constant bounded by 1.

We claim that $v$ is itself Lipschitz continuous in $x$ with the Lipschitz constant bounded by 1. To show this, take an $\varepsilon > 0$ and define $\tau_\varepsilon^\alpha$ as the first exit time of $(t, w_t + \alpha_t)$ from $Q_\varepsilon = (-\varepsilon, \varepsilon) \times B_\varepsilon$. Also, take a random variable $\xi$ which is uniformly distributed on $[0, 1]$ and independent of the filtration $\{\mathcal{F}_t\}$. Then for any $\delta \geq 0$, $\delta \xi$ is a stopping time with respect to the filtration $\{\mathcal{F}_t \vee \sigma(\xi)\}$. We also know that changing the probability space and filtrations does not affect the value function (see, e.g., [5, 6]). Set $\tau_{\varepsilon\delta}^\alpha = \tau_\varepsilon^\alpha \wedge (\delta\xi)$. If $(t, x)$ is such that $(t, x) + Q_\varepsilon \subset Q$, then $\tau_{\varepsilon\delta}^\alpha \leq \tau^{\alpha,t,x}$ for any control process $\alpha$. For $\alpha_t \equiv 0$, by Theorem 1.9,

$$(1.15) \qquad v(t, x) \geq R_{\varepsilon\delta} v(t, x) + S_{\varepsilon\delta} v(t, x) + P_{\varepsilon\delta} f(t, x),$$

where

$$R_{\varepsilon\delta} v(t, x) = E\{e^{-\delta\xi} v(t + \delta\xi, x + w_{\delta\xi}) I_{\delta\xi < \tau_\varepsilon^\alpha}\},$$

$$S_{\varepsilon\delta} v(t, x) = E\{e^{-\tau_\varepsilon^\alpha} v(t + \tau_\varepsilon^\alpha, x + w_{\tau_\varepsilon^\alpha}) I_{\delta\xi \geq \tau_\varepsilon^\alpha}\},$$

$$P_{\varepsilon\delta} f(t, x) = E \int_0^{\tau_{\varepsilon\delta}^\alpha} e^{-t} f(x + w_t)\, dt.$$

Obviously, as $\delta \downarrow 0$, $S_{\varepsilon\delta} v$ and $P_{\varepsilon\delta} f$ tend to zero uniformly with respect to $(t, x)$. Also, by the lower semicontinuity of $v$ and Fatou's lemma ($|v|$ is bounded),

$$\liminf_{\delta \downarrow 0} R_{\varepsilon\delta} v(t, x) \geq v(t, x) P(0 < \tau_\varepsilon^\alpha) = v(t, x),$$

which, along with (1.15), shows that $R_{\varepsilon\delta} v \to v$ as $\delta \downarrow 0$ on the set of $(t, x)$ such that $(t, x) + Q_\varepsilon \subset Q$. We now observe the obvious fact that the distribution of $(\delta\xi, w_{\delta\xi})$ has a density, so that in the definition of $R_{\varepsilon\delta}v$, we can replace $v$ with $\bar{v}$, which implies that

$$|R_{\varepsilon\delta}v(t,x) - R_{\varepsilon\delta}v(t,y)| \leq |x - y|$$

for all $\varepsilon > 0$ and $\delta > 0$. This and the above prove our claim.

Note that generally, since $g$ is only assumed to be lower semicontinuous, it is easy to see that $v$ need not be continuous in $\bar{D}$. Also, $v$ need not be continuous in $t$ unless $g$ is continuous in $t$.

EXAMPLE 1.12. In Example 1.11, the value function is most likely discontinuous because the data are unbounded. However, Assumption 1.7 does not guarantee that $v$ is continuous in $Q \cup \partial'Q$, even if everything is bounded and continuous. To see that, return to Example 1.4, letting (1.8) describe the dynamics for some control $\alpha \in A = \{1, 2\}$ and the equation

(1.16) $$dx_t = dw_t, \qquad dy_t = dz_t$$

describe the process response under the other control $\beta$, where $(w_t, z_t)$ is a two-dimensional Wiener process. Also, let $f^\alpha = f^\beta = 1$ and keep $g = 0$. Then, obviously, $v$ is greater than the function of the same name in Example 1.4. Also, $v(0, x, y) = 0$ for $(x, y) \in \partial B_2$ and, thus, $v(0, x, y)$ is discontinuous at the points $(-\sqrt{3}, \pm 1)$. However, as is easily checked the function

$$\psi(t, x, y) = 2(2 - r)(r - 1), \qquad r = \sqrt{x^2 + y^2},$$

satisfies Assumption 1.7 with $\underline{\alpha} = \beta$.

EXAMPLE 1.13. In Example 1.12, the function $v$ is discontinuous only at a few points on the boundary. One can modify this example in such a way that Assumption 1.7 is still satisfied and the discontinuities occur inside $Q$. To show that, replace (1.16) with

$$dx_t = b(r_t)x_t\, dt, \qquad dy_t = b(r_t)y_t\, dt,$$

where $b(r)$ is a smooth function on $[1, 2]$, with $b(r) = -1$ for $r \in [1, 5/4]$ and $b(r) = 1$ for $r \in [7/4, 2]$. Then any smooth function $\psi(t, x)$ such that $\psi = r - 1$ near $\partial B_1$ and $\psi = 2 - r$ near $\partial B_2$ satisfies Assumption 1.7 near $\partial B_1 \cup \partial B_2$. According to Remark 1.8, one can find a new function satisfying Assumption 1.7 as it is stated.

However, it is not hard to see that the new value function $v$ coincides with the function from Example 1.4 on the set where $t = 0$, $r \in [1, 5/4]$, $x < 0$ and $y \in [-1, 0]$. Since, on the other hand, $v$ is not less than the function from Example 1.4, $v(0, x, y)$ is discontinuous on that part of the line $(x, -1)$, $x \leq 0$, which lies in $B_{5/4} \setminus B_1$.

REMARK 1.14. Theorem 2.1 in Chapter V of [4] concerning Bellman's principle requires the existence of a rather smooth function $\bar{g}$ in $\bar{Q}$ such that $\bar{g} = g$ on $(-1, T) \times \partial D$, $\bar{g} \geq g$ on $\{T\} \times D$ and in $Q$,

$$(1.17) \qquad D_t \bar{g} + \sup_{\alpha \in A} [L^\alpha \bar{g} + f^\alpha] \leq 0.$$

This assumption is not satisfied in Examples 1.4 and 1.12 where $g \equiv 0$ because otherwise, by Itô's formula, we would have $v \leq \bar{g}$ in $Q \cup \partial' Q$ and $v(0, x, y)$ would go to zero as $(x, y)$ goes to $\partial D$.

**2. Lipschitz continuity of $v$ in $x$ and Hölder continuity of $v$ in $t$.** In this section, we show that under certain additional conditions, the function $v$ defined in Section 1 is Lipschitz continuous with respect to $x$ and Hölder $1/2$ continuous in $t$. Both the case that $Q$ is a general domain and the case that $Q$ is a cylindrical domain are treated.

For some applications (see, e.g., the proof of Theorem 2.17), it is also convenient to investigate the dependence of $v$ on parameters. Therefore, apart from our basic objects and assumptions introduced at the beginning of Section 1, we suppose that for an $\varepsilon_0 \in [0, 1]$ and each $\varepsilon \in \{0, \varepsilon_0\}$, we are also given

$$(2.1) \qquad \begin{aligned} \sigma^\alpha(\varepsilon) &= \sigma^\alpha(t, x, \varepsilon), & b^\alpha(\varepsilon) &= b^\alpha(t, x, \varepsilon), \\ c^\alpha(\varepsilon) &= c^\alpha(t, x, \varepsilon), & f^\alpha(\varepsilon) &= f^\alpha(t, x, \varepsilon), & g(\varepsilon) &= g(t, x, \varepsilon), \end{aligned}$$

having the same meaning and satisfying the same assumptions as the original $\sigma, b, c, f$. The solution of (1.1) corresponding to $\sigma^\alpha(\varepsilon), b^\alpha(\varepsilon)$ will be denoted by $x_s^{\alpha, t, x}(\varepsilon)$ and the functions $v^\alpha$, $v$ constructed from the new objects by $v^\alpha(t, x, \varepsilon)$ and $v(t, x, \varepsilon)$, respectively. We assume that for $\varepsilon = 0$, the functions in (2.1) coincide with the original ones, so that in our notation,

$$v^\alpha(t, x) = v^\alpha(t, x, 0), \qquad v(t, x) = v(t, x, 0).$$

Naturally, the operator $L^\alpha$ constructed from $\sigma^\alpha(\varepsilon), b^\alpha(\varepsilon)$ and $c^\alpha(\varepsilon)$ is denoted by $L^\alpha(\varepsilon)$, and by $\tau^{\alpha, t, x}(\varepsilon)$, we mean the first exit time of $(t + s, x_s^{\alpha, t, x}(\varepsilon))$, $s \geq 0$, from $Q$.

Let $\lambda \in [0, \infty), K, K_1, T \in (0, \infty)$ be constants. The names of the following assumptions contain a parameter $\varepsilon$. This is done in order to provide flexibility for using the assumptions in different settings.

ASSUMPTION 2.1 ($\varepsilon$). (i) We have $Q \subset (-\infty, T) \times \mathbb{R}^d$ and in $\bar{Q}$, we are given a continuous function $\psi$ such that $\psi = 0$ on the parabolic boundary $\partial' Q$ of $Q$.

(ii) The functions $g(\varepsilon)$ and $\psi$, their first and second derivatives in $x$ and first derivatives in $t$ are a continuous on $\bar{Q}$.

(iii) For each $\alpha \in A$ on $Q$, we have

$$|f^\alpha(\varepsilon) + D_t g(\varepsilon) + L^\alpha(\varepsilon) g(\varepsilon)| \leq K_1, \qquad c^\alpha(\varepsilon) \geq \lambda.$$

ASSUMPTION 2.2 ($\varepsilon$).  For any $\alpha \in A$, it holds that

$$D_t \psi + L^\alpha(\varepsilon) \psi \leq -1 \qquad \text{in } Q.$$

ASSUMPTION 2.3 ($\varepsilon$).  (i) For $\zeta(\varepsilon) = \sigma^\alpha(\varepsilon), b^\alpha(\varepsilon)$, $\alpha \in A$ and any $(t, x)$, $(t, y) \in Q$, we have

$$|\zeta(t, x, \varepsilon) - \zeta(t, y, 0)| \leq K(|x - y| + \varepsilon).$$

(ii) For $\zeta(\varepsilon) = \psi$, $c^\alpha(\varepsilon)$, $g(\varepsilon)$, $f^\alpha(\varepsilon)$, $\alpha \in A$ and any $(t, x), (t, y) \in Q$, we have

$$|\zeta(t, x, \varepsilon) - \zeta(t, y, 0)| \leq K_1(|x - y| + \varepsilon),$$

where $|\sigma|$ has the usual meaning (trace $\sigma\sigma^*$)$^{1/2}$ for matrices $\sigma$.

We start by estimating the moments of the difference of solutions of (1.1) with different initial values.

THEOREM 2.4.  *Let Assumption 2.3($\varepsilon$)(i) be satisfied for some $\varepsilon \in \{0, \varepsilon_0\}$. Take any $p \geq 0$, $(t, x), (t, y) \in Q$, $\alpha \in \mathfrak{A}$ and a stopping time $\gamma \leq \tau^{\alpha, t, x} \wedge \tau^{\alpha, t, y}(\varepsilon)$. Then*

$$(2.2) \qquad E \sup_{s \leq \gamma} e^{-Ms} |x_s^{\alpha, t, x} - x_s^{\alpha, t, y}(\varepsilon)|^p \leq 3(|y - x|^p + \varepsilon^p),$$

*where $M = M(p, K) \geq 0$.*

PROOF.  First, we take $p \geq 2$. For simplicity of notation, we drop the indices $\alpha, t, x, y$ in what follows. For instance, we denote $x_r = x_r^{\alpha, t, x}$. Also, set $y_r = x_r^{\alpha, t, y}(\varepsilon)$.

By using Itô's formula, for $s \in [0, \gamma]$, we obtain

$$e^{-Ms}(|x_s - y_s|^p + \varepsilon^p)$$
$$= |x - y|^p + \varepsilon^p + m_s$$
$$\quad - M \int_0^s (|x_r - y_r|^p + \varepsilon^p) e^{-Mr} \, dr$$
$$\quad + p \int_0^s |x_r - y_r|^{p-2}(x_r - y_r)^*(b(t + r, x_r) - b(t + r, y_r, \varepsilon)) e^{-Mr} \, dr$$
$$\quad + \frac{p(p-2)}{2} \int_0^s |x_r - y_r|^{p-2} |\sigma(t + r, x_r) - \sigma(t + r, y_r, \varepsilon)|^2 e^{-Mr} \, dr,$$

where $m_s$ is a local martingale starting at zero. Due to Assumption 2.3($\varepsilon$)(i), we can choose $M = M(p, K) \geq 1$ sufficiently large so that

$$e^{-Ms}(|x_s - y_s|^p + \varepsilon^p) \leq |x - y|^p + \varepsilon^p + m_s.$$

Upon applying Lemma 7.3(i) of [9], we get

$$Ee^{-M\gamma}(|x_\gamma - y_\gamma|^p + \varepsilon^p) \leq |x - y|^p + \varepsilon^p.$$

Since $\gamma$ is any stopping time $\leq \tau^{\alpha,t,x} \wedge \tau^{\alpha,t,y}(\varepsilon)$, by Lemma 7.3(ii) of [9], we conclude that

$$
\begin{aligned}
E \sup_{s \leq \gamma} e^{-M\delta s} |x_s - y_s|^{p\delta} &\leq E \sup_{s \leq \gamma} e^{-M\delta s}(|x_s - y_s|^p + \varepsilon^p)^\delta \\
&\leq \frac{2 - \delta}{1 - \delta}(|x - y|^p + \varepsilon^p)^\delta \\
&\leq \frac{2 - \delta}{1 - \delta}(|x - y|^{p\delta} + \varepsilon^{p\delta}) \\
&\leq 3(|x - y|^{p\delta} + \varepsilon^{p\delta})
\end{aligned}
$$

for any $\delta \in (0, 1/2)$. It only remains to observe that when $p$ runs through $[2, \infty)$ and $\delta$ through $(0, 1/2)$, the product $p\delta$ covers $(0, \infty)$.   $\square$

By using (1.11), we arrive at the following:

**Lemma 2.5.** *Let Assumptions 2.1($\varepsilon$) and 2.2($\varepsilon$) be satisfied for some $\varepsilon \in \{0, \varepsilon_0\}$. Then on $Q \cup \partial' Q$, we have*

$$|v(t, x, \varepsilon) - g(t, x, \varepsilon)| \leq K_1 \psi(t, x).$$

**Theorem 2.6.** *Let Assumptions 2.1(0) and 2.2(0) be satisfied. Take some $\varepsilon \in \{0, \varepsilon_0\}$ and suppose that Assumptions 2.1($\varepsilon$), 2.2($\varepsilon$) and 2.3($\varepsilon$) are satisfied. Then there are constants $N$ depending only on $K, K_1, d_1, d$ and $M$ depending only on $K$ such that for any $(t, x), (t, y) \in Q \cup \partial' Q$, we have*

$$(2.3) \qquad |v(t, x) - v(t, y, \varepsilon)| \leq N e^{(T-t)(M-\lambda)_+}(|x - y| + \varepsilon).$$

**Proof.** Due to Lemma 2.5, we may concentrate on points inside $Q$. Fix $(t, x), (t, y) \in Q$ and for any $\alpha \in \mathfrak{A}$ and $s \geq 0$, set

$$y_s^{\alpha,t,x} = x_s^{\alpha,t,y}(\varepsilon), \qquad \bar{\varphi}_s^{\alpha,t,x} = \int_0^s c^{\alpha_r}(t + r, y_r^{\alpha,t,x}, \varepsilon) \, dr.$$

This notation will allow us to use our convention regarding indices with which we provide the expectation sign.

By using Theorem 1.6 with

$$\gamma^\alpha = \tau^{\alpha,t,x} \wedge \tau^{\alpha,t,y}(\varepsilon),$$

we get

$$(2.4) \qquad |v(t,x) - v(t,y,\varepsilon)| \le I_1 + I_2,$$

where

$$I_1 = \sup_{\alpha \in \mathfrak{A}} E_{t,x}^{\alpha} \int_0^{\gamma} |f^{\alpha_s}(t+s,x_s) e^{-\varphi_s} - f^{\alpha_s}(t+s,y_s,\varepsilon) e^{-\bar{\varphi}_s}| \, ds,$$

$$I_2 = \sup_{\alpha \in \mathfrak{A}} E_{t,x}^{\alpha} |v(t+\gamma,x_\gamma) e^{-\varphi_\gamma} - v(t+\gamma,y_\gamma,\varepsilon) e^{-\bar{\varphi}_\gamma}|.$$

By using the inequality $|e^a - e^b| \le e^{a \vee b}|a-b|$ and Assumption 2.3, we obtain

$$|f^{\alpha_s}(t+s,x_s) e^{-\varphi_s} - f^{\alpha_s}(t+s,y_s,\varepsilon) e^{-\bar{\varphi}_s}|$$

$$\le N e^{-\lambda s} \left[ |x_s - y_s| + \varepsilon + \int_0^s (|x_r - y_r| + \varepsilon) \, dr \right]$$

$$\le N e^{(\mu-\lambda)s} \left[ e^{-\mu s}(|x_s - y_s| + \varepsilon) + \int_0^s e^{-\mu r}(|x_r - y_r| + \varepsilon) \, dr \right]$$

$$\le N e^{(\mu-\lambda)s}(1+s) \sup_{s \le \gamma} e^{-\mu s}(|x_s - y_s| + \varepsilon)$$

$$\le N e^{(\mu+1-\lambda)s} \sup_{s \le \gamma} e^{-\mu s}(|x_s - y_s| + \varepsilon),$$

where $\mu$ is any constant $\ge 0$. Upon applying Theorem 2.4, we get

$$(2.5) \qquad \begin{aligned} I_1 &\le N \sup_{\alpha \in \mathfrak{A}} E_{t,x}^{\alpha} \sup_{s \le \gamma} e^{-Ms}(|x_s - y_s| + \varepsilon) \int_0^{\gamma} e^{(M+1-\lambda)s} \, ds \\ &\le N(|x-y| + \varepsilon) e^{(M+2-\lambda)_+(T-t)}. \end{aligned}$$

To estimate $I_2$, we observe that either $(t+\gamma,x_\gamma)$ or $(t+\gamma,y_\gamma)$ is on $\partial' Q$. Due to Lemma 2.5, in the first case, we have

$$(2.6) \qquad \begin{aligned} &|v(t+\gamma,x_\gamma) e^{-\varphi_\gamma} - v(t+\gamma,y_\gamma,\varepsilon) e^{-\bar{\varphi}_\gamma}| \\ &= |g(t+\gamma,x_\gamma) e^{-\varphi_\gamma} - v(t+\gamma,y_\gamma,\varepsilon) e^{-\bar{\varphi}_\gamma}| \\ &\le |g(t+\gamma,x_\gamma) e^{-\varphi_\gamma} - g(t+\gamma,y_\gamma,\varepsilon) e^{-\bar{\varphi}_\gamma}| + N|\psi(t+\gamma,y_\gamma)| e^{-\gamma\lambda} \\ &\le N e^{(M+1-\lambda)\gamma} \sup_{s \le \gamma} e^{-Ms}(|x_s - y_s| + \varepsilon) \\ &\quad + N|\psi(t+\gamma,x_\gamma) - \psi(t+\gamma,y_\gamma)| e^{-\gamma\lambda} \\ &\le N e^{(M+1-\lambda)\gamma} \sup_{s \le \gamma} e^{-Ms}(|x_s - y_s| + \varepsilon). \end{aligned}$$

A similar argument is valid in the second case. Thus, by Theorem 2.4, for any $\alpha \in \mathfrak{A}$, we have

$$
(2.7) \qquad E_{t,x}^{\alpha} |v(t+\gamma, x_\gamma)e^{-\varphi_\gamma} - v(t+\gamma, y_\gamma, \varepsilon)e^{-\bar{\varphi}_\gamma}|
$$

$$
\leq N(|x-y| + \varepsilon)e^{(M+1-\lambda)_+(T-t)}.
$$

After combining (2.4), (2.5) and (2.7), we obtain (2.3) with $M+2$ in place of $M$. $\square$

THEOREM 2.7. *Under Assumptions 2.1(0), 2.2(0) and 2.3(0) also suppose that for any $\alpha \in A$ in $Q$,*

$$
(2.8) \qquad |\sigma^\alpha| + |b^\alpha| \le K_1.
$$

*Then there are constants $N = N(K, K_1, d, d_1)$ and $M = M(K)$ such that for any $(t, x), (s, x) \in Q \cup \partial' Q$ such that $|s - t| \le 1$, we have*

$$
(2.9) \qquad |v(s, x) - v(t, x)| \le |g(s, x) - g(t, x)| + K_1|\psi(s, x) - \psi(t, x)|
$$

$$
+ N|s - t|^{1/2} e^{(T-t)(M-\lambda)_+}.
$$

*In particular, if $g$ and $\psi$ are Hölder-1/2 continuous in $Q \cup \partial' Q$ with respect to $t$, then so is $v$.*

PROOF. Observe that if both points $(s, x)$ and $(t, x)$ are on $\partial' Q$, then the left-hand side of (2.9) is less than the first term on the right and so there is nothing to prove. However, if one of them is in $Q$, then $x$ can be slightly moved in such a way that they both fall into $Q$ and by Theorem 2.6, this leads to an insignificant modification of the left-hand side of (2.9). We see that it suffices to concentrate on $(s, x), (t, x) \in Q$.

Next, we assume that $t \le s$ and set $\gamma^{\alpha, t, x} = (s - t) \wedge \tau^{\alpha, t, x}$. Note that by Bellman's principle and Itô's formula (as usual, we drop indices $\alpha, t, x$ from objects behind the expectation sign),

$$
(v - g + K_1\psi)(t, x)
$$

$$
= \sup_{\alpha \in \mathfrak{A}} E_{t,x}^{\alpha} \Big[ (v - g + K_1\psi)(t+\gamma, x_\gamma)e^{-\varphi_\gamma}
$$

$$
+ \int_0^\gamma (\bar{f}^{\alpha_r} + K_1(D_t + L^{\alpha_r})\psi)(t+r, x_r)e^{-\varphi_r} \, dr \Big],
$$

where

$$
\bar{f}^\alpha = f^\alpha + (D_t + L^\alpha)g, \qquad |\bar{f}^\alpha| \le K_1, \qquad \bar{f}^\alpha + K_1(D_t + L^\alpha)\psi \le 0.
$$

Since $v - g + K_1\psi = 0$ on $\partial' Q$,

$$
E_{t,x}^{\alpha}(v - g + K_1\psi)(t+\gamma, x_\gamma)e^{-\varphi_\gamma} = E_{t,x}^{\alpha}(v - g + K_1\psi)(s, x_{s-t})e^{-\varphi_{s-t}} I_{\gamma = s-t},
$$

so that by Theorem 2.6,

$$E_{t,x}^{\alpha}(v - g + K_1\psi)(t + \gamma, x_\gamma)e^{-\varphi_\gamma}$$
$$\leq E_{t,x}^{\alpha}[(v - g + K_1\psi)(s,x) + Ne^{(T-t)(M-\lambda)_+}|x - x_{s-t}|]e^{-\varphi_{s-t}}I_{\gamma=s-t}.$$

Furthermore, well-known estimates of stochastic integrals, combined with the assumption that $\sigma$ and $b$ are bounded and that $s - t \leq \sqrt{s-t}$, imply that

$$E_{t,x}^{\alpha}|x - x_{s-t}| \leq N\sqrt{s-t}.$$

Next, according to Lemma 2.5, we have $(v - g + K_1\psi)(s,x) \geq 0$. It follows that

$$(v - g + K_1\psi)(t,x) \leq (v - g + K_1\psi)(s,x) + Ne^{(T-t)(M-\lambda)_+}\sqrt{s-t},$$
$$v(t,x) - v(s,x) \leq |g(t,x) - g(s,x)|$$
$$+ K_1|\psi(t,x) - \psi(s,x)| + Ne^{(T-t)(M-\lambda)_+}\sqrt{s-t}.$$

That

$$v(t,x) - v(s,x) \geq -|g(t,x) - g(s,x)|$$
$$- K_1|\psi(t,x) - \psi(s,x)| - Ne^{(T-t)(M-\lambda)_+}\sqrt{s-t}$$

is proved similarly by considering $v - g - K_1\psi$ and noting that this function is negative on $Q$. The theorem is thus proved.  □

Next, we consider the case where

$$Q = (-1, T) \times D$$

is a cylindrical domain in $\mathbb{R}^{d+1}$ under weaker assumptions on the boundary data. Let $D$ be a bounded domain and let $\psi(t,x)$, $g_1(\varepsilon) = g_1(t,x,\varepsilon)$ and $g_2(\varepsilon) = g_2(x,\varepsilon)$ be functions on $\bar{Q}$.

ASSUMPTION 2.8 ($\varepsilon$). (i) The functions $g_1(\varepsilon)$ and $\psi$, their first derivatives with respect to $(t,x)$ and their second derivatives with respect to $x$ are continuous on $\bar{Q}$, $\psi > 0$ in $Q$ and $\psi$ vanishes on $(-1,T) \times \partial D$.

(ii) We have

$$g(\varepsilon) = g_1(\varepsilon) \qquad \text{on } (-1,T) \times \partial D,$$

$$g(\varepsilon) = g_2(\varepsilon), \qquad |g_2(\varepsilon)| \leq K_1, \qquad |g_2(\varepsilon) - g_1(\varepsilon)| \leq K_1\psi \qquad \text{on } \{T\} \times \bar{D}.$$

(iii) For each $\alpha \in A$ on $Q$, we have

$$|f^\alpha(\varepsilon) + D_t g_1(\varepsilon) + L^\alpha(\varepsilon)g_1(\varepsilon)| \leq K_1, \qquad c^\alpha(\varepsilon) \geq \lambda.$$

Observe that

$$\partial' Q = ((-1, T) \times \partial D) \cup (\{T\} \times \bar{D}).$$

Itô's formula immediately yields the following:

LEMMA 2.9. *Let Assumptions* 2.8($\varepsilon$) *and* 2.2($\varepsilon$) *be satisfied for some* $\varepsilon \in \{0, \varepsilon_0\}$. *Then on* $Q \cup \partial' Q$, *we have*

$$(2.10) \qquad |v(\varepsilon) - g_1(\varepsilon)| \le K_1 \psi.$$

The following theorem can be proven in almost the same way as Theorem 2.6. By "the assertion of Theorem 2.6" in Theorem 2.10 we mean that which follows "Then" in the statement of Theorem 2.6. Theorem 2.13 should be read similarly.

THEOREM 2.10. *Let Assumptions* 2.8(0) *and* 2.2(0) *be satisfied. Take some* $\varepsilon \in \{0, \varepsilon_0\}$ *and suppose that Assumptions* 2.8($\varepsilon$), 2.2($\varepsilon$) *and* 2.3($\varepsilon$) *are satisfied if in Assumption* 2.3($\varepsilon$), *we replace* $g$ *with* $g_1, g_2$. *Then the assertion of Theorem* 2.6 *again holds true.*

Indeed, we can reproduce the proof of Theorem 2.6 except that we use Theorem 1.9 in place of Theorem 1.6 and while estimating $I_2$ instead of (2.6), we write

$$|v(t + \gamma, x_\gamma)e^{-\varphi_\gamma} - v(t + \gamma, y_\gamma, \varepsilon)e^{-\bar{\varphi}_\gamma}|$$
$$= \cdots I_{\gamma < T-t} + \cdots I_{\gamma = T-t}$$
$$= |g_1(t + \gamma, x_\gamma)e^{-\varphi_\gamma} - v(t + \gamma, y_\gamma, \varepsilon)e^{-\bar{\varphi}_\gamma}|I_{\gamma < T-t}$$
$$\quad + |g_2(x_{T-t})e^{-\varphi_{T-t}} - g_2(y_{T-t}, \varepsilon)e^{-\bar{\varphi}_{T-t}}|I_{\gamma = T-t},$$

where as before, the first term on the right is less than the last term in (2.6) and the second is majorated by $I_{\gamma = T-t}$ times

$$K_1(|x_{T-t} - y_{T-t}| + \varepsilon)e^{-\lambda(T-t)} + K_1^2 e^{-\lambda(T-t)} \int_0^{T-t} (|x_r - y_r| + \varepsilon)\,dr.$$

REMARK 2.11. In Theorem 2.10, we required $\psi$ to satisfy Assumption 2.2 in $Q$. As in Remark 1.8, one may show that we actually need this assumption only near $(-1, T) \times \partial D$.

Using Theorem 7.4 of [7] (or the corresponding results in [5]) and the above results immediately yield the following:

COROLLARY 2.12. *Suppose that the assumptions of Theorem 2.10 or Theorem 2.6 are satisfied with $\varepsilon = 0$. Then for any $\alpha \in A$, (1.7) holds true in $Q$ in the sense of generalized functions, that is, for any nonnegative $\chi \in C_0^\infty(Q)$,*

$$\int_Q v(-D_t - c^\alpha)\chi \, dt \, dx$$

$$+ \int_Q [\chi(b^\alpha)^i D_i v + \chi f^\alpha - (\chi D_i(a^\alpha)^{ij} + (a^\alpha)^{ij} D_i \chi) D_j v] \, dt \, dx \leq 0.$$

Our next result concerns the Hölder continuity of $v$ in $t$.

THEOREM 2.13. *Under Assumptions 2.8(0), 2.2(0) and 2.3(0), suppose that (2.8) holds for any $\alpha \in A$ in $Q$. Then the assertions of Theorem 2.7 are valid with $g_1$ in place of $g$ in (2.9) and $v$ is Hölder-1/2 continuous in $Q \cup \partial' Q$ with respect to $t$.*

The proof of this theorem follows that of Theorem 2.7 almost word for word; of course, we replace $g$ with $g_1$ in that proof.

In the following remark, we state an analog of one of the assertions of Theorem 2.1 in Chapter V of [4]. As everywhere in the article, we remain within the framework introduced in Section 1.

REMARK 2.14. Let $Q = (-1, T) \times D$ and let $\psi$ be a function on $\bar{Q}$ which is continuous along with its first derivatives with respect to $(t, x)$ and the second derivatives with respect to $x$. Also, let Assumption 2.2(0) be satisfied, let $\psi > 0$ in $Q$ and let $\psi = 0$ on $(-1, T) \times \partial D$. Assume that $A = A(1)$ and $g$ is continuous. It then turns out that $v$ is continuous in $\bar{Q} \setminus (\{-1\} \times \bar{D})$.

Indeed, the fact that $A = A(1)$ guarantees the validity of (2.8) and Assumption 2.3(0)(i). Furthermore, having in mind approximations using mollifiers, we may assume that $c$ and $f$ are Lipschitz continuous in $x$ uniformly with respect to other variables and that $g$ is infinitely differentiable (see more about this in [5]). Then it only remains to observe that $v$ is continuous in $\bar{Q} \setminus (\{-1\} \times \bar{D})$ by Theorems 2.10 and 2.13 [of course, in these theorems, we take $g_1(\varepsilon) = g(\varepsilon) = g$ and $g_2(\varepsilon) = g(T, \cdot)$].

Note that our requirement that Assumption 2.2 be satisfied is, in fact, very similar to condition (1.17) imposed in Theorem 2.1 in Chapter V of [4]. However, we only need it for $g \equiv 0$, albeit with 1 in place of $f^\alpha$.

Before stating out last result, the obtaining of which largely motivated this article, we take a $\delta \in (0, 1]$, define $B = \{x \in \mathbb{R}^d : |x| < 1\}$, $\Lambda = (-1, 0)$,

$$\hat{A} = A \times \Lambda \times B,$$

for $\beta = (\alpha, r, y) \in \hat{A}$, set

$$(\sigma^\beta, b^\beta, c^\beta, f^\beta)(t, x) = (\sigma^\alpha, b^\alpha, c^\alpha, f^\alpha)(t + \delta^2 r, x + \delta y),$$

and introduce the following:

ASSUMPTION 2.15.   (i)  Assumptions 2.8(0), 2.2(0) and 2.3(0) hold if we there replace $A, \sigma^\alpha, b^\alpha, c^\alpha, f^\alpha$ with $\hat{A}, \sigma^\beta, b^\beta, c^\beta, f^\beta$, respectively.

(ii) For each $\alpha \in A$ in $(-2, T) \times \mathbb{R}^d$, the functions $\sigma^\alpha$ and $b^\alpha$ are Lipschitz continuous in $x$ with constant $K$ and Hölder-1/2 continuous in $t$ with constant $K$, and the functions $c^\alpha \geq \lambda$ and $f^\alpha$ are Lipschitz continuous in $x$ with constant $K_1$ and Hölder-1/2 continuous in $t$ with constant $K_1$. Condition (2.8) is also satisfied in that domain.

(iii) The functions $\psi$ and $g_1$ are defined on $H := [0, T] \times \mathbb{R}^d$ and are Lipschitz continuous in $x$ with constant $K_1$ and Hölder-1/2 continuous in $t$ with constant $K_1$.

REMARK 2.16.   On account of Assumption 2.15(ii) and the boundedness of the derivatives of $\psi$ and $g_1$ entering the operators $L^\alpha$, obviously, Assumption 2.15(i) is satisfied for sufficiently small $\delta$ and somewhat modified $\psi$ and $K_1$ if Assumptions 2.8(0), 2.2(0) and 2.3(0) hold in their original form.

Recall that until now, continuity in $t$ has not been assumed for $\sigma, b, c, f$.

Introduce

$$H(\delta) = [0, T - \delta^2] \times \mathbb{R}^d,$$

$$D^\delta = \{x \in D : \operatorname{dist}(x, \partial D) > \delta\}, \qquad Q(\delta) = (0, T - \delta^2) \times D^\delta.$$

For multi-indices $\gamma = (\gamma_1, \ldots, \gamma_d)$, $\gamma_i = 0, 1, \ldots,$ as usual, we set

$$D_x^\gamma = D_1^{\gamma_1} \cdots D_d^{\gamma_d}, \qquad |\gamma| = \gamma_1 + \cdots + \gamma_d.$$

THEOREM 2.17.   Extend $v$ as $g_1$ in $H \setminus \bar{Q}$. Suppose that $Q(\delta) \neq \varnothing$. Then under Assumption 2.15, there exists an infinitely differentiable function $u^\delta$ defined on $H(\delta)$ such that for any $\alpha \in A$,

$$(2.11) \qquad D_t u^\delta + L^\alpha u^\delta + f^\alpha \leq 0$$

in $Q(\delta)$ and for any integers $m \geq 1, k, l \geq 0$ and multi-indices $\gamma$ such that $2k + l = m$ and $|\gamma| = l$,

$$(2.12) \qquad |D_t^k D_x^\gamma u^\delta(t, x)| \leq N e^{(T-t)(M-\lambda)_+} \delta^{1-m},$$

$$(2.13) \qquad |u^\delta(t, x) - v(t, x)| \leq N e^{(T-t)(M-\lambda)_+} \delta$$

in $H(\delta)$, where $N = N(K, K_1, d, d_1)$ and $M = M(K)$.

Proof. Let $v^\delta$ be constructed from $\hat{A}, \sigma^\beta, b^\beta, c^\beta, f^\beta, g$ in the same way as $v$ from $A, \sigma^\alpha, b^\alpha, c^\alpha, f^\alpha, g$. Extend $v^\delta$ as $g_1$ in $H \setminus \bar{Q}$. Note that by virtue of Assumption 2.15 and Theorem 2.10 (where we take $\varepsilon = \varepsilon_0 = \delta$), we have

$$(2.14) \qquad |v(t,x) - v^\delta(t,x)| \leq N e^{(T-t)(M-\lambda)_+} \delta$$

in $(-1, T] \times \mathbb{R}^d$. Furthermore, by Assumption 2.15 and Theorems 2.10 and 2.13,

$$(2.15) \qquad |v^\delta(t,x) - v^\delta(s,y)| \leq N e^{(T-t)(M-\lambda)_+}(|t-s|^{1/2} + |x-y|)$$

if $(t,x), (s,y) \in (-1, T] \times \mathbb{R}^d$ and $|t - s| \leq 1$.

Take a nonnegative function $\zeta \in C_0^\infty(\mathbb{R}^{d+1})$ with support in $\Lambda \times B$ and unit integral. Our goal is to prove that

$$u^\delta(t,x) := \delta^{-d-2} \int_{\mathbb{R}^{d+1}} v^\delta(s,y) \zeta(\delta^{-2}(t-s), \delta^{-1}(x-y)) \, ds \, dy$$

is a function we need.

Inequality (2.12) follows from (2.15) and elementary properties of mollifiers which also imply that

$$|u^\delta(t,x) - v^\delta(t,x)| \leq N e^{(T-t)(M-\lambda)_+} \delta$$

in $H(\delta)$. By recalling (2.14), we see that it only remains to prove (2.11).

By Corollary 2.12, for any nonnegative $\chi \in C_0^\infty(Q)$,

$$\int_Q v^\delta(t,x)(-D_t - c^\alpha(t + \delta^2 r, x + \delta y))\chi(t,x) \, dt \, dx$$

$$+ \int_Q [\chi(t,x)((b^\alpha)^i(t + \delta^2 r, x + \delta y) D_i v^\delta(t,x)$$

$$+ f^\alpha(t + \delta^2 r, x + \delta y))$$

$$- (\chi(t,x) D_i(a^\alpha)^{ij}(t + \delta^2 r, x + \delta y)$$

$$+ (a^\alpha)^{ij}(t + \delta^2 r, x + \delta y) D_i \chi(t,x)) D_j v^\delta(t,x)] \, dt \, dx \leq 0.$$

We here substitute $\chi(t + \delta^2 r, x + \delta y)$ in place of $\chi(t,x)$ and change variables $(t + \delta^2 r, x + \delta y) \to (t,x)$ to find that for any fixed $r \in (-1,0)$, $|y| \leq 1$ and nonnegative $\chi \in C_0^\infty(Q(\delta))$,

$$\int_{Q(\delta)} v^\delta(t - \delta^2 r, x - \delta y)(-D_t - c^\alpha)\chi(t,x) \, dt \, dx$$

$$+ \int_{Q(\delta)} [\chi(b^\alpha)^i(t,x) D_i v^\delta(t - \delta^2 r, x - \delta y) + \chi f^\alpha(t,x)$$

$(2.16)$

$$- (\chi D_i(a^\alpha)^{ij} + (a^\alpha)^{ij} D_i \chi)(t,x)$$

$$\times D_j v^\delta(t - \delta^2 r, x - \delta y)] \, dt \, dx \leq 0.$$

After multiplying (2.16) by $\zeta(r,y)$, integrating with respect to $(r,y)$ and using Fubini's theorem, we obtain

$$
\int_{Q(\delta)} u^\delta(-D_t - c^\alpha)\chi \, dt \, dx
$$

$$
(2.17) \qquad + \int_{Q(\delta)} [\chi(b^\alpha)^i D_i u^\delta + \chi f^\alpha - (\chi D_i(a^\alpha)^{ij} + (a^\alpha)^{ij} D_i \chi) D_j u^\delta] \, dt \, dx
$$

$$
\leq 0.
$$

Since $u^\delta$ is smooth and $\chi \in C_0^\infty(Q(\delta))$ is an arbitrary nonnegative function, after one more integration by parts, inequality (2.17) implies (2.11).

The theorem is thus proved. $\square$

**Acknowledgments.** The authors are sincerely grateful to the referees for very thorough reading and many suggestions helping improve the presentation.

SCHOOL OF MATHEMATICS
INSTITUTE FOR ADVANCED STUDY
I EINSTEIN DRIVE
PRINCETON, NEW JERSEY 08540
USA
E-MAIL: hjdong@ias.edu

DEPARTMENT OF MATHEMATICS
UNIVERSITY OF MINNESOTA
127 VINCENT HALL
MINNEAPOLIS, MINNESOTA 55455
USA
E-MAIL: krylov@math.umn.edu